\documentclass[11pt]{article}

\usepackage[cp1252]{inputenc} 
\usepackage[a4paper]{geometry}
\usepackage{lmodern}
\usepackage[T1]{fontenc}
\usepackage{amsmath,amssymb,amsthm}
\usepackage{graphicx}
\usepackage{textcomp}
\usepackage{enumerate}

\usepackage[usenames,dvipsnames]{pstricks}
\usepackage{epsfig}
\usepackage{pst-grad} 
\usepackage{pst-plot} 

\newtheorem{prop}{Proposition}[section]
\newtheorem{cor}[prop]{Corollary}
\newtheorem{lem}[prop]{Lemma}

\newtheorem{thm}{Theorem}[section]

\theoremstyle{definition}
\newtheorem{defi}[prop]{Definition}
\newtheorem{ex}[prop]{Example}

\newcommand{\R}{\mathbb{R}}

\newcommand{\Z}{\mathbb{Z}}
\newcommand{\N}{\mathbb{N}}

\newcommand{\Addresses}{{
  \bigskip
  \footnotesize

  M.~Bouljihad, \textsc{UMPA, ENS Lyon, 46 allée d'Italie, 69007 Lyon, France
   }\par\nopagebreak
  \textit{E-mail address}, M.~Bouljihad: \texttt{mohamed.bouljihad@ens-lyon.fr}

}}

\title{Rigidity for group actions on homogeneous spaces by affine transformations}

\author{Mohamed Bouljihad}

\date{}

\begin{document}
\maketitle

\begin{abstract}
We give a criterion for the rigidity of the action of a group of affine transformations of a homogeneous space of a real Lie group. Let $G$ be a real Lie group, $\Lambda$ a lattice in $G$, and $\Gamma$ a subgroup of the affine group Aff$(G)$ stabilizing $\Lambda$. Then the action of $\Gamma$ on $G/\Lambda$ has the rigidity property in the sense of S. Popa \cite{Pop}, if and only if the induced action of $\Gamma$ on $\mathbb{P}(\frak{g})$ admits no $\Gamma$-invariant probability measure, where $\frak{g}$ is the Lie algebra of $G$. This generalizes results of M. Burger \cite{Bur}, and A. Ioana and Y. Shalom \cite{Ioa}. As an application, we establish rigidity for the action of a class of groups acting by automorphisms on nilmanifolds associated to step 2 nilpotent Lie groups.
\end{abstract}

\section{Introduction}

\setlength{\parindent}{0cm}

Property (T) for locally compact groups was introduced by Kazhdan in 1967 \cite{Kaz}, and he used it  to prove that certain lattices are finitey generated. It was first defined in terms of unitary representations (see \cite{Kaz}, \cite{Bek}). Property (T) has since become a basic notion in various areas of mathematics such as group theory, ergodic theory and operator algebras. 

The first examples of groups with property (T) are the special linear groups SL$_n(\R)$, for $n\geq 3$. Lattices in SL$_n(\R)$ are examples of discrete groups with property (T), since this property is inherited by lattices. 

As for the relative property (T), it was first coined by Margulis for pairs of groups $(G,H)$ where $G$ is a locally compact group and $H$ a subgroup of $G$ \cite{Mar}. Yet, the relative property (T) of $(\R^2\rtimes$SL$_2(\R),\R^2)$  was implicitly used in Kazhdan's original paper.

Burger \cite{Bur} studied the relative property (T) for pairs of groups of the form $(A\rtimes\Gamma, A)$ where $A$ is a discrete abelian group, and $\Gamma$ is a subgroup of Aut$(A)$. He gave a sufficient condition on the dual action of $\Gamma$ on $\hat{A}$ to ensure that $(A\rtimes\Gamma, A)$ has the relative property (T). In the case where $A=\Z^n$ and $\Gamma$ is a  subgroup of SL$_n(\Z)$, it gives rise to an action on the $n$-torus $\mathbb{T}^n=\R^n/\Z^n$.  His result can be stated as follows.

\begin{thm}[\cite{Bur}]

The pair of groups $(\Z^n\rtimes\Gamma,\Z^n)$ has relative property (T) if and only if there is no $\Gamma$-invariant probability measure on the projective space $\mathbb{P}(\R^n).$ 
\end{thm} 

Actually, the fact that the condition above is necessary was noticed in \cite{Cor} (Proposition 3.1.9, which is primarily based on the so-called Furstenberg's lemma \cite{Fur}).


\smallskip
In the spirit of relative property (T) for pairs of groups, Popa \cite{Pop} introduced the notion of relative property (T) for inclusions $B\subset N$ of von Neumann algebra (also called rigidity of inclusion).

Given a free ergodic action of a countable group $\Gamma$ on a non-atomic standard probability space $(X,\mu)$ by measure preserving transformations, there is an associated canonical von Neumann algebra; the group measure space construction L$^\infty(X)\rtimes\Gamma$ (\cite{MvN}). Following \cite{Gab}, we say that the action $\Gamma\curvearrowright (X,\mu)$ has \emph{the property (T) relative to the space} if the inclusion L$^\infty(X)\subset$ L$^\infty(X)\rtimes\Gamma$ has the relative property (T). This notion first appeared in \cite{Pop}, where such actions are called rigid. For an ergodic theoretical equivalent formulation of this property, see Proposition $2.3.$ below. 

In the case of group actions by automorphims on abelian groups, this notion is linked to the Kazhdan-Margulis relative property (T) as follows.

\begin{prop}[\cite{Pop}, Proposition 5.1]
Let $A$ be a discrete abelian group and let $\Gamma$ be a subgroup of $\textrm{Aut}(A)$. Let $\mu$ be the Haar measure on $\hat{A}$, the dual space of $A$. Then the following assertions are equivalent : 

\begin{enumerate}
\item The action $\Gamma \curvearrowright (\hat{A},\mu)$ has the property (T) relative to the space;
\item The pair of groups $(A\rtimes\Gamma, A)$ has the relative property (T).
\end{enumerate}
\end{prop} 

The previous proposition, together with Theorem 1.1, shows that the action $\Gamma\curvearrowright (\mathbb{T}^n,\mu)$ of a subgroup of SL$_n(\Z)$ has the property (T) relative to the space if, and only if, there is no $\Gamma$-invariant probability measure on $\mathbb{P}(\R^n)$.

\smallskip

Popa (\cite{Pop}, Problem 5.10.2) raised the following natural question : characterize the countable discrete groups for which there exists free ergodic probability measure preserving actions on $(X,\mu)$ that have the property (T) relative to the space. This question is still open. Yet, some results were obtained. It has been showed in \cite{Gab}, Theorem 1.3, that every non-amenable free product of groups admits such actions. Besides, it was proved in \cite{Ioa2}, Theorem 4.3, that every non-amenable group admits free ergodic measure preserving actions with a weak form of the property (T) relative to the space.

\smallskip

The notion of property (T) relative to the space of group actions played a fundamental role in Popa's work on fundamental groups of von Neumann $II_1$ factors (\cite{Pop}). It was recently used as a key ingredient to produce actions of $\mathbb{F}_\infty$ (the free group on infinitely many generators) whose $II_1$ factors and orbit equivalence relations have prescribed fundamental group (\cite{PopVa}). 

Recall that probability measure preserving actions of two countable groups $\Gamma$ and $\Lambda$ on two standard probability spaces $(X,\mu)$ and $(Y,\nu)$ are orbit equivalent, if they define isomorphic partitions of the spaces into orbits. It is a standard result that free ergodic actions of infinite amenable groups are all orbit equivalent (\cite{OW}). As for non-amenable groups, it has been proven, using the property (T) relative to the space of actions of $\mathbb{F}_2$ on $\mathbb{T}^2$, that every countable non-amenable group has uncountably many non-orbit equivalent actions (\cite{GabPop}, \cite{Ioa2}, \cite{Eps}). Moreover, new examples of standard equivalence relations with trivial outer automorphisms groups were given thanks to this notion of property (T) relative to the space (see \cite{PopVa}, \cite{Gab}).

\smallskip

Ioana and Shalom studied in \cite{Ioa} the property (T) relative to the space for actions given by a group of translations on homogeneous spaces of algebraic groups. They gave mild conditions under which this action has the property (T) relative to the space. Here is one of their results.

\begin{thm}[\cite{Ioa}, Theorem D]

Let $G$ be a real algebraic group and $\Lambda\subset G$ a lattice. Let $\Gamma \subset G$ be a countable subgroup and denote by $H$ its Zariski closure. Assume that $H$ has no proper normal co-compact algebraic subgroup and no non-trivial homomorphism into $\R^*$. Let $\eta$ be  a probability measure on $G/\Lambda$ which is invariant under translations by elements from $\Gamma$.

If the centralizer of $\Gamma$ (equivalently, of $H$) in $G$ is finite, then the action of $\Gamma$ by translations on $(G/\Lambda,\eta)$ has the property (T) relative to the space.

In the case $\eta=m_{G/\Lambda}$, the converse is true : if the action $\Gamma \curvearrowright (G/\Lambda,\eta)$ has the property (T) relative to the space, then the centralizer of $\Gamma$ in $G$ is finite.
\end{thm}


In our main result below, we are able to generalize these results in three ways : 
\begin{itemize}
\item $G$ is an arbitrary real Lie group (and $\Lambda$ is a lattive in $G$),
\item $\Gamma$ is an arbitrary countable group of affine transformations of $G/\Lambda$ (see below for this notion),
\item we give a simple necessary and sufficient condition for the property (T) relative to the space of the action $\Gamma\curvearrowright (G/\Lambda,\mu)$.
\end{itemize}
Moreover, our result shows that Burger's neat characterization of rigidity of actions of groups automorphism on tori extends to this general setting.

To state our result, we need to introduce some notation. Let $G$ be a locally compact group and $\Lambda$ a lattice in $G$. Let Aut$(G)$ be the group of continuous automorphisms of $G$ and Aff$(G)=G\rtimes$Aut$(G)$ the group of affine transformations of $G$. We denote by Aut$_\Lambda(G)$ the subgroup of Aut$(G)$ of all $\sigma\in $Aut$(G)$ with $\sigma(\Lambda)=\Lambda$. The group of affine transformations of $G/\Lambda$ $$\mathrm{Aff}_\Lambda(G)=G\rtimes\mathrm{Aut}_\Lambda(G)$$ acts in a natural way on $G/\Lambda$. Let $\mu$ be the $G$-invariant probability measure on $G/\Lambda$. Then $\mu$ is fixed by every $\gamma\in$Aff$_\Lambda(G)$.

Assume now that $G$ is a real Lie group with Lie algebra $\frak{g}$. The group Aff$(G)$ acts linearly on $\frak{g}$. This gives rise to an induced action of Aff$(G)$ on the projective space $\mathbb{P}(\frak{g})$. Here is our main result.

\begin{thm}
Let $G$ be a real Lie group with Lie algebra $\frak{g}$ and let $\Lambda$ be a lattice in $G$. Let $G/\Lambda$ be the corresponding homogeneous space. Let $\mu$ be the $G$-invariant probability measure on $G/\Lambda$.
Let $\Gamma$ be a countable subgroup of $\mathrm{Aff}_\Lambda(G)$. 
Then the following properties are equivalent :

\begin{itemize}
\item[(i)] $\Gamma \curvearrowright (G/\Lambda,\mu)$ has the property (T) relative to the space;

\item[(ii)] there is no $\Gamma$-invariant measure probability on $\mathbb{P}(\frak{g})$;

\item[(iii)] the pair of groups $( \frak{g}^*\rtimes \Gamma,\frak{g}^*)$ has the relative property (T).

\end{itemize} 
\label{thm1}
\end{thm}

We should mention that our approach is largely based on techniques developed in \cite{Ioa}. 

Up until recently (see \cite{Ioa}), the examples of actions of group which have the property (T) relative to the space were built using the relative property (T) of a pair $(A\rtimes\Gamma, A)$. We give new examples of actions with the property (T) relative to the space which are not obtained from a pair of groups with relative property (T). For instance, let $N$ be a simply connected and connected free nilpotent Lie groups of step two. Let $\Lambda$ be a lattice of $N$. Assume that the dimension of $N$ is strictly bigger than $3$. Then (see Example \ref{ex1}) the action of Aut$_\Lambda(N)$ on $N/\Lambda$ has the property (T) relative to the space. 

\smallskip

The example given above is part of family of nilmanifolds associated to finite graphs which are studied in \cite{Dan}. We use Theorem \ref{thm1} to give a fairly simple criterion for actions of groups of automorphims of these nilmanifolds.



Let $(S,E)$ be a finite non-oriented graph. Let $N$ be the associated step two nilpotent Lie group and $N/N_\Z$ the corresponding standard nilmanifold. Let $V$ be a a complementary subspace  to $[N,N]$ in $N$. The group Aut$(N)$ has a decomposition as semi-direct product $T\rtimes G$ where $T$ (the group of the so-called shear automorphisms of $N$) is  isomorphic to some $\R^n$ and $G$ is a certain algebraic subgroup of GL(V) (see [DM05]). Let $L$ be the semi-simple Levi component of $G$. Every element $n\in N$ can be written $n=(x,y)$ where $x\in V$ and $y\in [N,N]$. Then $(t,u)\in T\rtimes L$ acts on $N$ as follows, $$(t,u).(x,y)=(u(x),P(u)(y)+t(u(x))),$$ where $P(u)\in$ GL$([N,N])$ is a certain matrix whose coefficients are polynomials in those of $u$ (see Section 4.2 for details).


In this setting, we have the following result.

\begin{thm} 

Let $\Gamma$ be a countable subgroup of $T\rtimes L$ stabilizing $N_\Z$. Assume $\Gamma$ is Zariski dense in $T\rtimes L$. The following assertions are equivalent

\begin{enumerate}[(i)]
\item $\Gamma\curvearrowright (N/N_\Z,\mu)$ has the property (T) relative to the space, where $\mu$ is the $N$-invariant probability measure on $N/N_\Z$,
\item the linear action of $\Gamma$ on $[\mathcal{N},\mathcal{N}]$ has no non-zero fixed vector, where $\mathcal{N}$ denotes the Lie algebra of $N$.
\end{enumerate}
\label{thm2}
\end{thm}

From this result, we deduce the following surprising fact. Let $H_\R$ be the Heisenberg group, $H_\Z$ the standard lattice. Despite the similarities with the action of the torus, the action $\text{Aut}_{H_\Z}(H_\R) \curvearrowright H_\R/H_\Z$ doesn't have the property (T) relative to the space (see Corollary \ref{cor}). In particular, any subequivalence relation of the equivalence relation induced by the action of $\text{Aut}_{H_\Z}(H_\R)$ on $H_\R/H_\Z$ doesn't have the property (T) relative to the space. 

\medskip

The paper is organized as follows. Section 2 contains a few preliminary results. The proof of Theorem \ref{thm1} is given in Section 3. In Section 4, we discuss the nilmanifols associated to finite graphs and prove Theorem \ref{thm2}.

\smallskip

\textbf{Acknowledgments.} 
I owe many thanks to Bachir Bekka for suggesting this problem, and for his helpful advices and comments. I am very grateful to Damien Gaboriau for his useful comments. I would also like to thank Mickael de la Salle for pointing out a mistake in a preliminary version of this paper.

\medskip
\section{Some preliminary lemmas}

The following lemma which appears in \cite{Ioa} will be useful in the sequel.

\begin{lem}[\cite{Ioa}]
Let $(X,\mu)$ a standard probability space \\
Let $c>0$ and $\eta_n$ a sequence of probability measure on $X\times X$ with the following properties:
\begin{itemize}
 \item $p_{*}^{1}(\eta_n) \leq c \mu $, where $p^i :X\times X\rightarrow X$ is the projection on the $i$-th factor ($i=1,2$), \\
 \item $\eta_n(\Delta_X)= 0 $, for all $x\in \N$,\\
 \item $\eta_n (A\times (X  \backslash  A ) ) \underset{n\rightarrow\infty}{\longrightarrow} 0$, for every Borel subset $A$ of  $X$.\\
\end{itemize}
Then, for every Borel partition $(A_i)_{i=1}^{\infty}$ of X, we have $$\eta_n(\bigcup_{i=1}^{\infty}(A_i\times A_i) )\,  \underset{n\rightarrow\infty}{\longrightarrow} \, 1.$$
\end{lem}

Our second lemma is an extension of Lemma F from \cite{Ioa}. The setting for this result and for the rest of this section is as follows. Let $G$ be a real Lie group, $\frak{g}$ be the Lie algebra of $G$ and denote by $p : \frak{g}  \backslash \{ 0\} \rightarrow \mathbb{P}(\frak{g})$ the canonical projection onto the projective space $\mathbb{P}(\frak{g})$ of $\frak{g}$. 

Let $q: G \rightarrow \frak{g}$ be any Borel map which equals the logarithm in some neighborhood $U$ of $e \in G$. Define $r : (G\times G)  \backslash  \Delta \rightarrow G,\,(x,y)\mapsto xy^{-1}$. 

As in $\cite{Ioa}$,we consider the map $$\rho\, : \,(G\times G)  \backslash  \Delta \rightarrow \mathbb{P}(\frak{g}),\,\rho(x,y)=p(q(r(x,y))).$$ We will denote by Ad$: G \rightarrow GL(\frak{g})$ the adjoint representation of $G$ on $\frak{g}$ as well as the associated action $G\rightarrow $Aut$(\mathbb{P}(\frak{g}))$ of $G$ on $\mathbb{P}(\frak{g})$.

Recall that the affine group Aff$(G)$ is the semi-direct product $G\rtimes$ Aut$(G)$. 

Now, we adapt the Lemma F \cite{Ioa} to our setting.

\begin{lem}

Let $X$ a Borel subset of $G$, and let $\mu$ be a probability measure on $X$. \\
Let $c>0$, and let $\eta_n$ be a sequence of probability measures on $X\times X$ with the following properties :
\begin{itemize}
 \item $p_{*}^{1}(\eta_n) \leq c \mu $, 
 \item $\eta_n(\Delta_X)= 0 $, for all $n \geq 1$, 
 \item $\eta_n (A\times (X  \backslash  A ) ) \underset{n\rightarrow\infty}{\longrightarrow} 0$, for every Borel subset A of X.
\end{itemize}

Let $\Lambda$ be a countable subgroup of $G$,  $\varphi : X \rightarrow \Lambda$ a Borel map and $\gamma=(\sigma,g)\in $Aff$(G)$ with $\sigma\in\mathrm{Aut}(G)$ and $g\in G$.\\
Set 
\[ D_{\gamma} = \{(x,y)\in X\times X  \backslash  \Delta_X  \,|\, \rho (\gamma(x)\varphi(x),\gamma(y)\varphi(y)) = \mathrm{Ad}(g)\circ\mathrm{d}\sigma_e(\rho (x,y)) \}\]
where $\mathrm{d}\sigma_e$ is the derivated automorphism at the identity of $\sigma$. Then
\[\lim_{n\rightarrow\infty}\eta_n (D_{\gamma}) =1.\]

\end{lem}

\bigskip

\textbf{Proof.} By definition of $q$ and of the adjoint representation Ad of $G$, we can find a neighborhood $V$ of the group unit $e$ in $G$ such that

\[q(g \sigma(x) g^{-1})=\mathrm{Ad}(g)\circ \mathrm{d}\sigma_e (q(x)),\quad\text{for all } x\in V.\] 

Set $A=\{(x,y)\in X\times X  \,|\,  xy^{-1} \in V \}$. \\

Choose a neighbourhood $W$ of $e$ in $G$ such that $WW^{-1} \subset V$. Since $G$ is separable, there exists a sequence $(h_i)_{i\in \mathbb{N}}$ of elements of $G$ such that $$G=\cup_{i=1}^{\infty} W h_i.$$
For every $i\geq 1$, let $$A_i= (W h_i  \backslash  (\bigcup_{j=1}^{i-1} W h_j)) \cap X.$$

Then $(A_i)_{i=1}^{\infty}$ is a Borel partition of $X$ and, since $\bigcup _{i=1}^{\infty }(A_i\times A_i) \subset A$, the previous lemma implies that $\eta_n(A) \underset{n\rightarrow\infty}{\longrightarrow}1$.

\medskip 

Let $B=\{ (x,y) \in X\times X \,|\, \varphi(x)=\varphi(y)\}$.

Since we can write $X$ as a countable partition $$X=\bigcup_{\lambda \in \Lambda} \{x  \,|\,  \varphi(x) = \lambda \}, $$
we have again by the previous lemma $$\eta_n(B)  \underset{n\rightarrow\infty}{\longrightarrow} 1.$$

\medskip

We claim that $A\cap B \subset D_{\gamma} \cup \Delta_X. $.

Indeed,  let $(x,y)\in A\cap B$ with $x\neq y$. Since $(x,y)\in B$, we have by definition of $\rho$
\[\begin{split}
 \rho(\gamma(x)\varphi(x),\gamma(y)\varphi(y)) &= \rho (\gamma(x)\varphi(x),\gamma(y)\varphi(x)) \\
  &= p(q(\gamma(x) \varphi(x) \varphi(x) ^{-1} \gamma(y) ^{-1} )) \\
  &= p(q(\gamma (xy^{-1}))) \\
  &= p(q(g \sigma(xy^{-1})g^{-1})). \\
\end{split}
\]  
Since $(x,y)\in A$, we have $xy^{-1}\in V$ and it follows that  
\[\begin{split}  
\rho(\gamma(x)\varphi(x),\gamma(y)\varphi(y))   &= p(\mathrm{Ad}(g)\mathrm{d}\sigma_e ( q((xy^{-1})))\\
  &= \mathrm{Ad}(g)\circ\mathrm{d}\sigma_e(\rho (x,y)).
  \end{split}
\]
Since $ \eta_n (A \cap B) \underset{n\rightarrow\infty}{\longrightarrow} 1$ and since $\eta_n(\Delta_X)=0$, we obtain that $ \eta_n(D_{\gamma}) \underset{n\rightarrow\infty}{\longrightarrow}1.$This finishes the proof of the lemma.

\begin{flushright}
 $\square$
\end{flushright}

\smallskip

We are now almost ready to prove the main result. We recall that Ioana (see \cite{Ioa3} Theorem 4.4 ,\cite{Ioa} Prop. 1) gave the following purely ergodic theoretical characterization of the property (T) relative to the space. It appears to be, in our case, much more handful than the original characterization given in term of pair of von Neumann algebra.

\begin{prop}[\cite{Ioa}]
\label{prop}
A measure preserving action $\Gamma\curvearrowright (X,\mu)$ of a countable group $\Gamma$ on a probability space $(X,\mu)$ has the property (T) relative to the space if and only if, for any sequence of Borel probability measures $\nu_n$ on $X\times X$ satisfying :
\begin{enumerate}
\item $p_*^i\nu_n=\mu$ for all $n$ and $i=1,2$, where $p^i:X\times X\rightarrow X$ denotes the projection onto the i-th coordinate,
\item $\int_{X\times X} \phi(x)\psi(y)d\nu_n(x,y)\underset{n\rightarrow\infty}{\longrightarrow}\int_X \phi(x)\psi(x)d\mu(x)$, for all bounded Borel functions $\phi,\psi$ on $X$,
\item $||(\gamma\times \gamma)_*\nu_n-\nu_n||\underset{n\rightarrow\infty}{\longrightarrow}0$, for every $\gamma\in\Gamma$,
\end{enumerate} we have that $\nu_n(\Delta_X) \underset{n\rightarrow\infty}{\longrightarrow} 1$ where $\Delta_X$ denotes the diagonal in $X\times X$.
\end{prop}

Note that in view of Proposition 1.1, a similar characterization was independently obtained in \cite{Tes}. The authors actually studied the relative property (T) for pair of groups of the form $(H\ltimes A, A)$, where $A$ is a locally compact abelian group, and $H$ is a locally compact group acting on $A$.

\section{Proof of Theorem \ref{thm1}}

The equivalence between the two last assertions follows from \cite{Cor}, Prop 3.1.9. Hence, we are only interested in the first two points. 

\smallskip

We proceed by contraposition : we will show that the action $\Gamma\curvearrowright (G/\Lambda,\mu) $ doesn't have the property (T) relative to the space if and only if there exists a $\Gamma$-invariant probability on $\mathbb{P}(\frak{g})$.

\medskip

\textbf{Direct implication :} Assume that $\Gamma \curvearrowright (G/\Lambda,\mu)$ doesn't have the property (T) relative to the space.

\medskip

Choose a Borel fundamental domain $X\subset G$ for the action of $\Lambda$ on G by right translations. 

As mentionned in the introduction, Aff$_\Lambda(G)$ acts in a natural way on $G/\Lambda$. We transfer this to an action of Aff$_\Lambda(G)$ on $X$ by setting

\[ \gamma . x = g\sigma(x)\omega (\gamma,x) \quad\text{for }x\in X,\,\gamma=(g,\sigma)\in\text{Aff}(G/\Lambda),\] where $\omega(\gamma,x)$ is the unique element in $\Lambda$ such that $g\sigma(x)\omega(\gamma,x)\in X$. 
 We denote by B$(X)$ the algebra of bounded measurable complex valued functions on $X$.

\medskip
Since $\Gamma\curvearrowright (X,\mu)$ doesn't the property (T) relative to the space, there exists a sequence of probability measures $\nu_n$ on $X\times X$ with the following properties 

\begin{enumerate}
\item $p_*^i\nu_n=\mu$, for all $n$ and $i=1,2,$
\item $\int_{X\times X}f(x)g(y)\mathrm{d}\nu_n(x,y)\underset{n\rightarrow\infty}{\longrightarrow} \int_X f(x)g(x)\mathrm{d}\mu$, for all $f,g\in\mathrm{B}(X),$
\item $\parallel (\gamma\times\gamma)_* \nu_n -\nu_n \parallel \underset{n\rightarrow\infty}{\longrightarrow} 0$, for all $\gamma\in \Gamma,$
\item $\underset{n\rightarrow \infty}{\liminf}\, \nu_n(\Delta_X) < 1.$ 
\end{enumerate}

Upon passing to a subsequence, we can assume that $c=\inf_{n\in\N} c_n > 0$, where $c_n=1-\nu_n(\Delta_X)$.

Define a sequence of probability measures $\eta_n$ on $X\times X$ by 
\[\eta_n(A)=\frac1{c_n}\nu_n(A \backslash \Delta_X),\quad\text{ for all Borel subsets } A\subset X\times X .\]
Properties 1 and 2 above imply that
\[p_*^1\eta_n \leq \frac1{c}\mu,\quad\text{ for all }n,\] 
and $\eta_n(A\times(X \backslash  A))\rightarrow 0$ for all Borel subset $A\subset X$.

Let $\gamma\in\Gamma$. Set $\varphi  \equiv \omega(\gamma,. )$ where $\omega(\gamma,.):X\rightarrow\Lambda,\,x\mapsto \omega(\gamma,x)$ is the Borel cocycle defined above. Lemma 2.2 shows that

\[\eta_n(\{(x,y)\in(X\times X) \backslash \Delta_X  \,|\,   \rho(\gamma . x, \gamma . y) =\gamma.(\rho(x,y))\} )\underset{n\rightarrow\infty}{\longrightarrow} 1.\quad (*)\]On the other hand, Property 3 implies that 
\[\parallel (\gamma\times\gamma)_*\eta_n-\eta_n\parallel\underset{n\rightarrow\infty}{\longrightarrow}0 .\]

\medskip

Consider the sequence of probability measure on $\mathbb{P}(\frak{g})$ defined by $\zeta_n=\rho_*\eta_n$. It follows from $(*)$ that we again have
\[\parallel \gamma_*\zeta_n-\zeta_n\parallel\underset{n\rightarrow\infty}{\longrightarrow} 0.\quad (**) \]

Now, since $\mathbb{P}(\frak{g})$ is a compact metric space, the space $\mathcal{M}(\mathbb{P}(\frak{g}))$ of probability measures on $\mathbb{P}(\frak{g})$, endowed with the weak-$*$ topology, is a compact metrizable space. So, passing to a subsequence, we can assume that $\zeta_n$ converges to a probability measure $\zeta$ on $\mathbb{P}(\frak{g})$. 

By $(**)$, $\zeta$ is $\gamma$-invariant. Since this holds for every $\gamma\in\Gamma$, we have found a $\Gamma$-invariant probability measure on $\mathbb{P}(\frak{g})$.

\bigskip

\textbf{Converse implication :} Assume now that there exists a $\Gamma$-invariant probability measure $\zeta$ on $\mathbb{P}(\frak{g})$.

\medskip We claim that the action of $\Gamma$ on $G/\Lambda$ doesn't have the property (T) relative to the space.

Denote by $H$ the Zariski closure of $\Gamma$ in GL$(\frak{g})$. By Theorem $3.11$ in \cite{Sha}, there exists a normal cocompact algebraic subgroup $H_0$ of the algebraic group $H$ such that $H_0$ fixes every point in the support supp$(\zeta)$ of $\zeta$.

Since supp$(\zeta)\neq \emptyset$, there exists $Y\in\frak{g} \backslash \{0\}$ and a group homomorphism $\chi:H_0\rightarrow\mathbb{R}^*$ such that
\[h.Y=\chi(h)Y,\quad \text{for all }h\in H_0.\]

It follows that $Y$ is invariant under the closure $\overline{[H_0,H_0]}$ of the commutator subgroup  of $H_0$. 
Consider $x_n=\exp(Y/n)$ for all $n\in\N$. Then for all $n\in \N$, $x_n\in G$ is fixed by $\overline{[H_0,H_0]}$ for the action given by :

\[\mathrm{A}(h_0)(x_n)=\exp(h_0.Y/n)=x_n, \,\forall n\in \N,\,\forall h_0\in[H_0,H_0]. \]

In this setting, we need to introduce the following action of Aff$(G)$ on $G$ by :

\[\mathrm{A} : \mathrm{Aff}(G)\rightarrow \mathrm{Aut}(G)\] with
\[\mathrm{A}(h)(x) = g\sigma(x) g^{-1}\,\text{ for } h=(\sigma,g)\in\text{Aff}(G),\, x\in G. \]

So that \[\mathrm{A}(h)(x_n)=\exp(h.Y/n),\text{ where } h\in \text{Aff}(G),\, \text{for all } n\in\N.\]

Since $x_n\rightarrow e$ and $x_n\neq e$ for $n$ large enough, and since $\Lambda$ is discrete, we can assume that $x_n\notin \Lambda$ for every $n$.

Consider the map
\[
\begin{split}
\varphi_n:G/\Lambda & \rightarrow G/\Lambda \times G/\Lambda \\
x & \mapsto (x , x_n x)
\end{split}
\]
for every $n$ and set $\nu_n = \varphi_n^*\mu$. Then $\nu_n$ is a probability measure on $G/\Lambda\times G/\Lambda$, and for all $\varphi \in B(G/\Lambda\times G/\Lambda)$ we have

\[\nu_n(\varphi)=\int_{x\in G/\Lambda} \varphi(x,x_n x)\text{d}\mu(x).\]

The measure $\nu_n$ being not necessarily invariant under $\Gamma$, we define a new sequence of measures obtained by averaging $\nu_n$ as follows.

First notice that $H/\overline{[H_0,H_0]}$ is amenable since $(H/\overline{[H_0,H_0]})/ (H_0/\overline{[H_0,H_0]})\cong H/H_0$ is amenable and $H_0/\overline{[H_0,H_0]}$ is an abelian normal subgroup of $H/\overline{[H_0,H_0]}$. So, by the Reiter condition (see \cite{Bek}, Appendix G), there exists $f_n\in L^1_{+,1}(H/\overline{[H_0,H_0]},\lambda)$ such that  $||\pi(h).f_n -f_n||\rightarrow 1$ for all $h\in H$, where $\pi$ is the regular representation of $H$ on L$^2(H/[H_0,H_0])$ and $\lambda$ is a Haar measure on $H/\overline{[H_0,H_0]}$.

We introduce the sequence of measures $\eta_{n,m}$ on $G/\Lambda\times G/\Lambda$ defined for $\varphi\in B(G/\Lambda\times G/\Lambda)$ by 

\[\eta_{n,m}(\varphi)=\int_{x\in G/\Lambda}\int_{h\in H/\overline{[H_0,H_0]}}f_m(h) \varphi(x,\mathrm{A}(h)(x_n) x)\text{d}\lambda(h)\text{d}\mu(x).\]

We claim that there exists a subsequence $(\eta_m)_m$ of $(\eta_{n,m})_{n,m}$ which satisfies properties 1.-3. of Proposition $\ref{prop}$, and such that $\underset{m\rightarrow \infty}{\liminf}\, \eta_{m} (\Delta_{G/\Lambda} )< 1$. 

\medskip

The projections of $\eta_{n,m}$ on each coordinate are equal to $\mu$, for all $n,m$, so that property 1 is verified for every subsequence. 

\medskip

To check property 3, let $\gamma \in \Gamma$. Write $\gamma=(\sigma,g)$ where $\sigma\in$ Aut$(G)$ and $g\in G$. For all $\varphi\in B(G/\Lambda\times G/\Lambda)$, we have :

\[
\begin{split}
(\gamma\times \gamma)_*\eta_{n,m}(\varphi) =&\int_{x\in G/\Lambda}\int_{h\in H/\overline{[H_0,H_0]}} f_m(h)\varphi(\gamma(x),\gamma(\mathrm{A}(h)(x_n) x))\text{d}\lambda(h)\text{d}\mu(x) \\
 =&\int_{x\in G/\Lambda}\int_{h\in H/\overline{[H_0,H_0]}} f_m(h)\varphi(\gamma(x),g\sigma((\mathrm{A}(h)(x_n))\sigma(x))\text{d}\lambda(h)\text{d}\mu(x),\end{split}\]
 
 since $\sigma$ is an automorphism. Therefore,
 
 \[\begin{split}
 (\gamma\times \gamma)_*\eta_{n,m}(\varphi)=&\int_{x\in G/\Lambda}\int_{h\in H/\overline{[H_0,H_0]}} f_m(h)\varphi(\gamma(x),g\sigma (\mathrm{A}(h)(x_n)) g^{-1}g\sigma(x))\text{d}\lambda(h)\text{d}\mu(x)\\
 =&\int_{x\in G/\Lambda}\int_{h\in H/\overline{[H_0,H_0]}} f_m(h)\varphi(\gamma(x),g\sigma (\exp(h.Y/n)) g^{-1}\gamma(x))\text{d}\lambda(h)\text{d}\mu(x)\\
 =& \int_{x\in G/\Lambda}\int_{h\in H/\overline{[H_0,H_0]}} f_m(h) \varphi(x,\exp((\gamma h). Y/n) x)\text{d}\lambda(h)\text{d}\mu(x),\end{split}\]

since  $\mu$ is $\Gamma$-invariant.

As such,
 \[\begin{split}
 (\gamma\times \gamma)_*\eta_{n,m}(\varphi) =& \int_{x\in G/\Lambda}\int_{h\in H/\overline{[H_0,H_0]}} f_m(\gamma^{-1}.h') \varphi(x,\exp(h'. Y/n) x)\text{d}\lambda(h)\text{d}\mu(x),\\ 
 =& \int_{x\in G/\Lambda}\int_{h'\in H/\overline{[H_0,H_0]}} f_m(\gamma^{-1}.h')) \varphi(x,\mathrm{A}(h') (x_n) x)\text{d}\lambda(h')\text{d}\mu(x),\end{split}\]
 
 by the $H$-invariance of $\lambda$. Hence,
 
\[\begin{split}
| (\gamma\times \gamma)_*\eta_{n,m}(\varphi)-\eta_{n,m}(\varphi) | \leq& |\int_{x\in G/\Lambda} \int_{h\in H/\overline{[H_0,H_0]}} (f_m(\gamma^{-1}.h) - f_m(\gamma) )\varphi(x,\mathrm{A}(h)(x_n)x) \mathrm{d}\lambda(h) \mathrm{d}\mu(x)| \\
\leq& \int_{h\in H/\overline{[H_0,H_0]}} |f_m(\gamma^{-1}.h)-f_m(h)| \big(\int_{x\in G/\Lambda}|\varphi(x,\mathrm{A}(h)(x_n)x)|\mathrm{d}\mu(x)\big)\mathrm{d}\lambda(h) \\
\leq & || \pi(\gamma).f_m-f_m||\end{split}\]

In summary, we obtain

\[ || (\gamma\times \gamma)_*\eta_{n,m}-\eta_{n,m}|| \leq || \pi(\gamma).f_m-f_m|| .\] Therefore, property 3. will be verified as soon as $m$ tends to infinity.

\medskip

We now verify property 2. It is enough to check this property for functions of the form $\varphi\otimes\psi$, where $\varphi$ and $\psi$ are continuous functions on the standard Borel space $G/\Lambda$.

For $\varphi$, $\psi\in$ C$(G/\Lambda)$, we check using Lebesgue's dominated convergence theorem that, for any fixed $m\in\N$, we have $$\int_{h\in H/\overline{[H_0,H_0]}}|f_m(h)|\int_{x\in G/\Lambda}|\psi(\mathrm{A}(h)(x_n)x)-\psi(x)|\mathrm{d}\mu(x)\text{d}\lambda(h)\rightarrow 0,$$ since $x_n\rightarrow e$ when $n\rightarrow +\infty$. Therefore,

\[
\begin{split}
\int_{(x,y)\in G/\Lambda\times G/\Lambda}\varphi(x)\psi(y)\mathrm{d}\eta_{n,m}(x,y)&=\int_{x\in G/\Lambda}\int_{h\in H/\overline{[H_0,H_0]}}f_m(h)\varphi(x)\psi(\mathrm{A}(h)(x_n)x)\text{d}\lambda(h)\mathrm{d}\mu(x)\\
&\underset{n\rightarrow\infty}{\longrightarrow}\int_{x\in G/\Lambda}\varphi(x)\psi(x)\mathrm{d}\mu(x).\end{split}\] 
We claim that there exists a subsequence $\eta_m$ that satisfies $\underset{m\rightarrow \infty}{\liminf}\, \eta_{m} (\Delta_{G/\Lambda} )< 1$. Let $X$ be a fundamental domain for the action of $G$ on $G/\Lambda$. We can assume that $\mu_G(X) =1$, where $\mu_G$ is the Haar measure on $G$. Let $U$ be a fixed neighborhood of $e$ in $G$ such that $\Lambda\cap U =\{e\}$. Then we have for any fixed $m\in \N$, using once again Lebesgue's dominated convergence theorem

\[ \int_{h\in H/\overline{[H_0,H_0]}} f_m(h) \mu_G(\{x\in X  \,|\,  x^{-1}\mathrm{A}(h)(x_n)x \notin U \})\mathrm{d}\lambda(h)\rightarrow 0  \] as $n\rightarrow \infty$.

Let $(\varphi_n\otimes\psi_n)_n$ a sequence of dense functions in C$(G/\Lambda)\otimes $C$(G/\Lambda)$. For every $m\in \N$, there exists $n_1(m)\in\N$ such that for all $n\geq n_1(m)$ we have

\[ |\int_{G/\Lambda\times G/\Lambda}\varphi_1(x)\psi_1(y)\mathrm{d}\eta_{n,m}(x,y)-\int_{G/\Lambda}\varphi_1(x)\psi_1(x)\mathrm{d}\mu(x)|\leq \frac1{2^m}\] and \[\int_{h\in H/\overline{[H_0,H_0]}} f_m(h) \mu(\{x\in X  \,|\,  x^{-1}\mathrm{A}(h)(x_n)x \notin U \})\mathrm{d}\lambda(h)\leq \frac1{2^m}.\]
In the same way, for every $m\in \N$, there exists $n_2(m)\geq n_1(m)$ such that for all $n\geq n_2(n_1(m))$ we have

\[ |\int_{G/\Lambda\times G/\Lambda}\varphi_2(x)\psi_2(y)\mathrm{d}\eta_{n,m}(x,y)-\int_{G/\Lambda}\varphi_2(x)\psi_2(x)\mathrm{d}\mu(x)|\leq \frac1{2^m}.\]

We finally consider the subsequence of $\eta_{n,m}$ defined for all $m\in \N$ by $\eta_m=\eta_{n_m(m),m}.$

One can easily check that this subsequence satisfies properties $1.-3.$.

\medskip

We claim that $\eta_m(\Delta_{G/\Lambda})$ does not converge to $1$. Indeed, assume by contradiction that $\eta_m(\Delta_{G/\Lambda})\rightarrow 1$, i.e. $$\int_{h\in H/\overline{[H_0,H_0]}} f_m(h) \mu(\{x\in G/\Lambda  \,|\,  x=\mathrm{A}(h)(x_{n(m)})x \})\mathrm{d}\lambda(h)\underset{m\rightarrow\infty}{\longrightarrow} 1.$$

In particular for $X$ the fundamental domain for the action of $G$ on $G/\Lambda$ fixed before, we have 

 $$\int_{h\in H/\overline{[H_0,H_0]}} f_m(h) \mu_G(\{x\in X  \,|\,  x^{-1}\mathrm{A}(h)(x_{n(m)})x\in \Lambda \})\mathrm{d}\lambda(h)\underset{m\rightarrow\infty}{\longrightarrow}  1.$$

Recall that by construction of our subsequence we have 

\[ \int_{h\in H/\overline{[H_0,H_0]}} f_m(h) \mu_G(\{x\in X  \,|\,  x^{-1}\mathrm{A}(h)(x_{n(m)})x \notin U \})\mathrm{d}\lambda(h)\underset{m\rightarrow\infty}{\longrightarrow} 0, \]where $U$ is a neighborhood of $e$ such that $\Lambda \cap U=\{e\}$.  

Thus,

\[ \int_{h\in H/\overline{[H_0,H_0]}} f_m(h) \mu_G(\{x\in X  \,|\,  x^{-1}\mathrm{A}(h)(x_{n(m)})x \in \Lambda\cap U \})\mathrm{d}\lambda(h)\underset{m\rightarrow\infty}{\longrightarrow}  1  \]

i.e.\[ \mu_G(\{ x\in X  \,|\,  x_{n(m)} = e \}) \underset{m\rightarrow\infty}{\longrightarrow} 1\] 

This implies that $x_{n(m)}=e$ for $m$ large enough. This contradicts the fact that $x_{n(m)}=\exp(\frac{Y}{n(m)})\neq e$ for large $m$.

Hence $\eta_m(\Delta_{G/\Lambda})$ doest not converge to 1.

This shows that the action of $\Gamma$ on $G/\Lambda$ doesn't have the property (T) relative to the space and the proof is complete.
\begin{flushright}
$\square$
\end{flushright}

\section{Actions on nilmanifolds associated with graphs}

We give new examples of actions which have the property (T) relative to the space which do not arise from pairs of groups with relative property (T). Our examples are based on a class of nilpotent Lie groups associated with graphs which appear in \cite{Dan}.

\subsection{Action on the Heisenberg nilmanifold}

First, let us consider the realization of the Heisenberg group $H_\R$ as $\R^3$ with the law 

\[ (x,y,z).(x',y',z')=(x+x',y+y',z+z'+(xy'-x'y)) ,\quad \forall (x,y,z),(x',y',z')\in \R^3 \]

The subset $H_\Z=\Z^3$ is a lattice. We consider the associated nilmanifold $H_\R / H_\Z$ with the Lebesgue measure.

The group of automorphism of $H_\R$ stabilizing $H_\Z$ is (see \cite{Heu},) 
\[\text{Aut}_{H_\Z}(H_\R) =T\rtimes L\]

where 

\[
\begin{array}{llll} 
T&=\{\begin{pmatrix} I_2 &  0 \\ t & 1  \end{pmatrix}  \,|\,  t\in\text{Hom}(\Z^2,\Z)\}, &\quad L&=\{ \begin{pmatrix} g & 0 \\0  & 1 \end{pmatrix}  \,|\,  g\in\text{SL}_2(\Z)\} \\
&\cong \text{Hom}(\Z^2,\Z) & &\cong \text{SL}_2(\Z) 
\end{array}
\]

The action of $(t,g)\in\text{Aut}(H_\R) $, with $t\in\text{Hom}(\Z^2,\Z), g\in \text{SL}_2(\Z)$ on $H_\R$ is given by  

\[ (t,g).(x,y,z)=(g(x,y),z+t(x,y)), \quad \text{for } (x,y,z)\in H_\R \]

Observe that every $g\in \text{Aut}_{H_\Z}(H_\R)$ fixes pointwise the image of the center of $H_\R$ in $H_\R/H_\Z$. From Theorem \ref{thm1}, we deduce the following corollary 

\begin{cor}
For every subgroup $\Gamma$ of $\text{Aff}_{H_\Z}(H_\R)$, the action $\Gamma \curvearrowright H_\R/H_\Z$ doesn't have the property (T) relative to the space.
\label{cor}
\end{cor}

\subsection{Nilpotent Lie groups associated with graphs}

Let $(S,E)$ be a finite non oriented graph, with $S$  the set of vertices and $A$ the set of edges. The edge between the vertices $\alpha$ and $\beta$ will be denoted by $\alpha \beta$.

Let $V$ be the $\R$-vector space with $S$ as a basis, and $W$ be the $\R$-vector space with $E$ as a basis. Consider $\mathcal{N}=V\oplus W$. A Lie algebra structure is defined on $\mathcal{N}$ by

\[ \begin{split}
[u_\alpha,u_\beta]&= u_{\alpha\beta}\in W\,\text{ if } \, \alpha,\beta\in S\\
&= 0\, \text{ otherwise }
\end{split}\] where $u_\alpha,u_\beta $ denote basis vectors in $V$ and $W$.

It is clear that $\mathcal{N}$ is a 2-step nilpotent Lie algebra.

The simply connected and connected Lie group $N$ with Lie algebra $\mathcal{N}$ can be realized as $\mathcal{N}$ with the multiplication defined by 

\[ (v_1,w_1).(v_2,w_2)=(v_1+v_2,w_1+w_2+\frac1{2}[v_1,v_2])\] for all $v_1,v_2\in V$ and $w_1,w_2\in W$.

Observe that the exponential map is simply the identity.

\begin{ex} 
\hfill

\begin{enumerate}[(i)]
\item The Heisenberg group is given by the graph with two vertices and one edge.
\item Free nilpotent Lie groups of step 2 are given by complete graphs.
\end{enumerate}
\end{ex}

Let $N_\Z$ be the subgroup of $N$ generated by $\{(v,w)\in N  \,|\,  v\in \Z^{|S|},\,w\in\Z^{|E|}    \}$. Then $N_\Z$ is a cocompact lattice in $N$; we call $N / N_\Z$ the compact nilmanifold associated with $(S,E)$. 

Observe that Aut$(N)$=Aut$(\mathcal{N})$.

We choose a complementary subspace $V$ of $[\mathcal{N},\mathcal{N}]$ so that $\mathcal{N}=V\oplus [\mathcal{N},\mathcal{N}]$ and denote by $\pi : \mathcal{N} \rightarrow V$ the associated projection. We define a group monomorphism

\[\begin{split}
\text{Hom}(V,[\mathcal{N},\mathcal{N}]) & \rightarrow \text{Aut}(\mathcal{N}) \\
\theta &\mapsto t_\theta 
\end{split}
\] by $t_\theta(x)=x +\theta(\pi(x))$, $x\in \mathcal{N}$. Set

\[T= \{ t_\theta  \,|\,  \theta \in \text{Hom}(V,[\mathcal{N},\mathcal{N}]) \} \]
and
\[ G = \{ g\in \text{Aut}(\mathcal{N})  \,|\, g(V)=V\} .\]

It is straightforward to check that 

\[ \text{Aut}(\mathcal{N})= T\rtimes G\]  (see \cite{Dan}).

With respect to the decomposition $\mathcal{N}=V\oplus [\mathcal{N},\mathcal{N}]$, we can write 

\[\begin{split}
T&= \{\begin{pmatrix} I & 0 \\ t & I \end{pmatrix}  \,|\, t\in\text{Hom}(V,[\mathcal{N},\mathcal{N}])\} \\
G &=\{ \begin{pmatrix} g & 0 \\ 0 & P(g) \end{pmatrix}  \,|\,  g\in\text{GL}(V))\}
\end{split}\] where for $g\in $ GL$(V)$, $P(g)\in$ GL$([\mathcal{N},\mathcal{N}])$ is a certain matrix whose coefficients are polynomials in those of $g$.

As a result, to understand Aut$(\mathcal{N})$ one has to study $G$, and for this, we recall the following definition from \cite{Dan}.

\begin{defi}[\cite{Dan}]
A subset $S'\subset S$ is said \emph{coherent} if for every pair $\alpha,\beta\in S'$, and for every $\gamma \in S'$ such that $\alpha\gamma\in E$, either $\gamma =\beta$ or $\beta\gamma \in E$.
\end{defi}

We can decompose $S$ into \emph{maximal coherent components}, defining this way an equivalence relation on $S$. We denote by $[S]$ the set of coherent components, and for any $[\alpha]\in [S]$, let $V_{[\alpha]}$ be the subspace of $V$ generated by the elements of $[\alpha]$. We can then write $V=\bigoplus_{[\alpha]} V_{[\alpha]}$.

The following description of $G$ was obtained in \cite{Dan} : the connected component of the identity of $G$ is 

\[(\prod_{[\alpha]\in[S]} GL^+(V_{[\alpha]})). M\] for a closed nilpotent connected normal subgroup $M$ of $G$.

In the sequel, we will only consider subgroups of $G$ contained in 
\[ L\cong\prod_{[\alpha]\in[S]} SL(V_{[\alpha]}) \]

Let us look at a few examples.

\begin{ex}
\hfill

\begin{enumerate}[(i)]

\item The Heisenberg group arises from the graph with 2 vertices and 1 edge ; it has only one coherent component. So, in this case, $L=$SL$_2(\mathbb{R})$ and

\[T\rtimes L \cong \text{Hom}(\R^2,\R)\rtimes \text{SL}_2(\R)
\]

\item The complete graph with $n\geq 3$ vertices has as automorphism group

\[T\rtimes L \cong \text{Hom}(\R^n,\R^{\frac{n(n-1)}{2}})\rtimes \text{SL}_n(\R) 
\]

\end{enumerate}
\end{ex}

\subsection{Proof of Theorem \ref{thm2}}

Recall that we have the following notations. Let $\mathcal{N}$ be the nilpotent Lie algebra associated to a finite graph $(S,E)$. Let $N$ be the corresponding simply connected nilpotent Lie group and $N/N_\Z$ the associated manifold.

Write $\mathcal{N}=V\oplus [\mathcal{N},\mathcal{N}]$ for $V$ is complementary subspace to $[\mathcal{N},\mathcal{N}]$ in $\mathcal{N}$ and set $T\cong$ Hom$(V,[\mathcal{N},\mathcal{N}])$ and $L = \prod_{[\alpha]\in[S]} SL(V_{[\alpha]}) $ as in section 4.2.

Observe that $L=\prod_{[\alpha]\in[S]} SL(V_{[\alpha]})$ has no non-trivial homomorphism into $\R^*$ and has no non-trivial cocompact normal subgroup.

Let us recall as well how $T\rtimes L$ acts on $N$ identified with $\mathcal{N}$ :

\[(t,g).(Y_1,Y_2)=(g(Y_1), P(g)(Y_2)+t(g(Y_1))) \quad (*)\] for $(t,g)\in T\rtimes L$ and $(Y_1,Y_2)\in \mathcal{N} = V\oplus [\mathcal{N},\mathcal{N}]$

\medskip

We prove the theorem by contraposition : we will show that $\Gamma$ has a non-zero fixed vector $Y\in[\mathcal{N},\mathcal{N}]$ if and only if the action of $\Gamma$ on $N/N_\Z$ doesn't have the property (T) relative to the space.

\medskip

\textbf{Direct implication :} Assume that $\Gamma$ has a non-zero fixed vector $Y\in[\mathcal{N},\mathcal{N}]$. 

The Dirac measure corresponding to the image of $Y$ in $\mathbb{P}(\mathcal{N})$ is then fixed by $\Gamma$. Hence, $\Gamma\curvearrowright N/N_\Z$ doesn't have the property (T) relative to the space by Theorem \ref{thm1}.

\medskip

\textbf{Converse implication :} Assume that the action of $\Gamma$ on $N/N_\Z$ doesn't have the property (T) relative to the space. By Theorem \ref{thm1}, $\Gamma$ fixes a probability measure $\xi$ on $\mathbb{P}(\mathcal{N})$. The stabilizer of $\xi$ in $\overline{\Gamma}^{\text{Zar}}=T\rtimes L$ is an algebraic group, and contains a normal cocompact subgroup $H$ which fixes every point in supp$(\xi)$ (\cite{Sha}, Theorem 3.11).

Choose $Y\in \mathcal{N}$ such that its image in $\mathbb{P}(\mathcal{N})$ belongs to supp$(\xi)$, and write $Y=(Y_1,Y_2)$ with $Y_1\in V$, $Y_2\in [\mathcal{N},\mathcal{N}]$. Then

\[(t,g).Y=\lambda(t,g) Y,\quad \forall (t,g)\in H\quad (**)\] for a homomoprhism $\lambda:H\rightarrow\mathbb{R}^*$. We claim that $Y_1=0$. 

\medskip
Indeed, recall that, with respect to the decomposition $V=\bigoplus_{[\alpha]\in[S]} V_{[\alpha]}$, we have

\[L = \begin{pmatrix} SL(V_{[\alpha_1]}) &  &  &  & &\\  & \ddots &  & & &\\ & & SL(V_{[\alpha_k]}) & & & \\ & & & I_{[\alpha_k]} & & \\& & & & \ddots & \\ & & & & & I_{[\alpha_r]} \end{pmatrix}.\]

Write $Y_1=(Y_{1,1},\dots,Y_{1,r})$ with $Y_{1,i}$ the component of $Y_1$ in $V_{[\alpha_i]}$.
 Assume first that $Y_{1,i}\neq 0$ for some $1\leq i \leq k$. Since SL$(V_{[\alpha_i]})$ is simple, and $H\cap$SL$(V_{[\alpha_i]})\neq \{1\}$, we have SL$(V_{[\alpha_i]})\subset H$; this is a contradiction as SL$_n(\mathbb{R})$ stabilizes no line in $\R^n$. So $$Y_{1,1}=\dots=Y_{1,k}=0.$$ It remains to show that $$Y_{1,k+1}=\dots=Y_{1,r}=0.$$

Assume by contradiction that $Y_{1,i}\neq 0$ for some $k+1 \leq i \leq r$. Then, looking at the $i$-th component in equation $(**)$, we see that $\lambda : H\rightarrow \mathbb{R}^*$ is trivial. It follows from $(*)$ and $(**)$ that 

\[ Y_2 = P(g)Y_2 + tY_1,\quad \forall (t,g)\in T\rtimes L.\]

This is a polynomial equation in the matrix coefficients of $t$ and $g$; since it is verified for $H$, it holds for its Zariski closure denoted by $M$ :

\[ Y_2 = P(g)Y_2+tY_2,\quad \forall (t,g)\in M \quad (***).\]

Observe that $M$ is still a normal cocompact subgroup of $T\rtimes L$. We claim that $M=T\rtimes L$. Once proved, it will follow from $(***)$ applied to $g=e$ that 

\[ Y_2=Y_2+tY_1,\quad \forall t\in T\]
and hence $Y_1=0$.

To show that $M=T\rtimes L$, observe that the projection of $M$ on $L$ is $L$, since it has no non-trivial normal cocompact subgroup. So, it suffices to prove that $T\subset M$.

For this, consider the group isomorphism $\varphi : T/ T\cap M \rightarrow T M / M$. Since $T M=T\rtimes L$, and since $T/T\cap M $ and $(T\rtimes L) /M$ are $\sigma$-compact locally compact groups, the group isomorphism $\varphi$, which is continuous, is a homeomorphism.

The group $T/T\cap M\cong T M/M$ is compact (as $M$ is cocompact in $T M$); so $T\cap M$ is a normal cocompact Zariski closed subgroup of $T=$Hom$(V,[\mathcal{N},\mathcal{N}])$ which is isomorphic to some $\R^m$. Hence, $T\subset M$. So, we have proved that $Y_1=0$. 

We claim now that $Y_2$ is fixed by $\Gamma$. Since $Y_2\in [\mathcal{N},\mathcal{N}]$, this will finish the proof. To show this, let $p_2 :T\rtimes L\rightarrow L$ denote the canonical projection. Recall that $L$ has no non-trivial homomorphism into $\R^*$ and has no non-trivial cocompact normal subgroup. It follows that $p_2(H)=L$. This shows that $Y_2$ is fixed by $L$ and hence by $\Gamma$.

\begin{flushright}
$\square$
\end{flushright}

\medskip

We will denote by $(T\rtimes L)_\Z$ the subgroup of $T\rtimes L$ whose elements have integers matrix coefficients.

\medskip

\begin{ex}

\begin{enumerate}[(i)]

\item Let $N$ be the nilpotent Lie group associated with the complete graph with $n\geq 3$ vertices. Consider the following subgroup of automorphisms

\[T\rtimes L \cong  \begin{pmatrix} I_n & 0  \\ \text{Hom}(\R^n,\R^{\frac{n(n-1)}{2}}) & I_{\frac{n(n-1)}{2}}  \end{pmatrix} \rtimes \begin{pmatrix} \text{SL}_n(\R) & 0 \\ 0  & P(\text{SL}_n(\R)) \end{pmatrix}\]

where $P$ is an injective morphism. The action of $(T\rtimes L)_\Z$ on $\mathcal{N}$ has no fixed point. Hence, the action of $(T\rtimes L)_\Z$ on $N/N_\Z$ have the property (T) relative to the space, by Theorem \ref{thm2}.

\item Let $N$ be the nilpotent Lie group associated with the "star" graph; this is the graph given by a central vertex which is linked with $n$ other vertices. Consider the following subgroup 

\[T\rtimes L \cong  \begin{pmatrix} I_{n+1} & 0  \\ \text{Hom}(\R^{n+1},\R^n) & I_{n}  \end{pmatrix} \rtimes \begin{pmatrix} \text{SL}_n(\R) &  0 & 0 \\0 &1 & 0\\ 0&0 &\text{SL}_n(\R)\end{pmatrix}\] of Aut$(N)$, where SL$_n(\R)$ acts on $[\mathcal{N},\mathcal{N}]$. There is obviously no fixed point in $[\mathcal{N},\mathcal{N}]$. Hence, the action of $(T\rtimes L)_\Z$ on $N/N_\Z$ has the property (T) relative to the space, by Theorem \ref{thm2}.

\end{enumerate}

\label{ex1}

\end{ex}

\bigskip

\nocite{Bek}
\nocite{Ioa}
\nocite{Pop}
\nocite{Sha}
\nocite{Con}
\nocite{Bur}
\nocite{Dan}
\nocite{Cor}
\nocite{Ioa2}
\nocite{Ioa3}
\nocite{Gab}
\nocite{GabPop}
\nocite{PopVa}
\nocite{Eym}
\nocite{Bek1}
\nocite{Fur}
\nocite{Dan1}
\nocite{Ioa4}
\nocite{Heu}
\nocite{Kaz}
\nocite{Tes}

\bibliographystyle{alpha}

\begin{small}
\bibliography{/Users/Moms/Dropbox/LaTeX/Bibliographie/bib}
\end{small}

\Addresses

\end{document}